\documentclass[11pt,
]{article}
\usepackage{amssymb,amstext,amsthm,latexsym}
\usepackage{mathrsfs,mathbbol}

\usepackage{amsmath}
\usepackage{amsfonts,amssymb}

\usepackage{mparhack}

\newcommand{\nek}{\newcommand}
\nek{\renek}{\renewcommand}

\DeclareMathAlphabet{\cur}{U}{eur}{m}{n}
\nek{\skr}{\mathscr}

\nek{\bfsf} {\sffamily\bfseries\upshape}

\nek{\ui}{\fontshape{ui}\selectfont}
\nek{\bfit}{\bfseries\itshape}
\nek{\bfsl}{\bfseries\slshape}
\nek{\sfbs}{\mdseries\sffamily\itshape}        %\bfseries

\nek{\vyk} [1] {}

\nek{\imar}[1]{\marginpar[%\vspace{-1ex}%
\flushright\footnotesize%
$\mtho\longrightarrow$\\%
\vspace{-1ex}{#1}\vspace*{1ex}]%
{%\vspace{-1ex}%
\flushleft\footnotesize%
$\mtho\longleftarrow$\\%
\vspace{-1ex}{#1}\vspace*{1ex}}%
}%
\nek{\imae}[1]{\marginpar[%\vspace{-1ex}%
\flushright\footnotesize\vspace{-4ex}%
$\mtho\longrightarrow$\\%
\vspace{-1ex}{#1}$\mtho$\vspace*{2ex}]%
{%\vspace{-1ex}%
\flushleft\footnotesize\vspace{-4ex}%
$\mtho\longleftarrow$\\%
\vspace{-1ex}{#1}$\mtho$\vspace*{2ex}}
}%

\nek{\parf}{\setcounter{subsection}0%\setcounter{equation}0
\chapter}
\nek{\punk} {\subsection}
\renek{\thesubsection}
{\arabic{subsection}}

\nek{\ppunk} [1]{\subsubsection{#1}}
\renek{\thesubsubsection}
{\thesubsection\hspace*{0.1ex}\Asbuk{subsubsection}}
\vyk{
\makeatletter
\renek\subsection{\@startsection{subsection}{2}{\z@}%
                                     {3.25ex\@plus 1ex \@minus .2ex}%
                                     {1.5ex \@plus .2ex}%
                                     {\normalfont\large\bfseries\S\hspace{0.2ex}}}
\makeatother
}

\vyk{
\chardef\eurl ='025
\chardef\eurM ='026
\nek{\fla} {{\cur{\eurl}}} 
\nek{\fmu} {{\cur{\eurM}}} 
\chardef\eurF = '011
\nek{\fFi} {{\cur{\eurF}}}
}

\newcounter{enuf}
\nek{\enufi}{\addtocounter{enuf}{1}}
\nek{\fenu}{
\def\theenumi{$\mtho(\fnsymbol{enuf})$}

\enufi\itsep
}
\newcounter{enuc}
\nek{\enuci}{\addtocounter{enuc}{1}}
\nek{\cenu}{
\def\theenumi{$\mtho\arabic{enuc}^\circ$}

\enuci\itsep
}

\newcounter{enuF}
\nek{\enuFi}{\addtocounter{enuF}{1}}
\nek{\Fenu}{
\def\theenumi{\psur(\arabic{enuF})\psur}

\nek{\ifla}{\addtocounter{enuF}1\itla}
\enuFi\itsep
}

\nek{\itsep}{\itemsep=0.3ex plus 0.1ex minus 0.1ex}
\nek{\tenu}[1]{
\def\theenumi{#1}
\itsep
}
\nek{\tenui}[1]{
\def\theenumii{#1}
\itsep
}

\theoremstyle{plain}

\newtheorem{theore}             {Theorem}  
\newtheorem{corollar}  [theore]{Corollary}
\newtheorem{propo}     [theore]{Proposition}
\newtheorem{lemm}      [theore]{Lemma}

\theoremstyle{definition}
\newtheorem{defn}      [theore]{Definition}
\newtheorem{rem}       [theore]{Remark}
\newtheorem{que}       [theore]{Problem}
\newtheorem{prim}       [theore]{Example}
\newtheorem{bla}       [theore]{Blanket Agreement}

\newtheorem*{prF}{Proof}               %\def\theprF{}
\newtheorem*{aq} {Acknowledgements}          %\def\theaq{}

\nek{\thsp}{\hspace{0.1ex plus \mathsurround}}
\nek{\bpro}{\begin{propo}}
\nek{\epro}{\end{propo}}
\nek{\bcl} {\begin{cla}}
\nek{\ecl} {\end{cla}}
\nek{\bup}{\begin{upr}}
\nek{\eup}{\qed\end{upr}}
\nek{\euP}{\end{upr}}
\nek{\bct} {\begin{clt}}
\nek{\ect} {\end{clt}}
\nek{\bco} {\begin{con}}
\nek{\eco} {\qeDD{Construction}\end{con}}
\nek{\bcor}{\begin{corollar}}
\nek{\ecor}{\end{corollar}}
\nek{\baq} {\begin{aq}}
\nek{\eaq} {\end{aq}}
\nek{\bex} {\begin{prim}}
\nek{\eex} {\qeD\end{prim}}
\nek{\eeX} {\end{prim}}
\nek{\bdf} {\begin{defn}} %\rm\thsp}
\nek{\eDf} {\end{defn}}
\nek{\edf} {\qeD\end{defn}}
\nek{\edF} {\end{defn}}
\nek{\ble} {\begin{lemm}}
\nek{\ele} {\end{lemm}}
\nek{\bsle} {\begin{subl}}
\nek{\esle} {\end{subl}}
\nek{\bte} {\begin{theore}}
\nek{\ete} {\end{theore}}
\nek{\bpra} {\begin{prav}} 
\nek{\epra} {\qeD\end{prav}} 
\nek{\bre} {\begin{rem}} 
\nek{\ere} {\qeD\end{rem}} 
\nek{\eRe} {\end{rem}} 
\nek{\bqe} {\begin{que}} 
\nek{\eqe} {\qeD\end{que}} 
\nek{\bql} {\begin{quel}} 
\nek{\eql} {\qeD\end{quel}} 
\nek{\bqo}{\begin{quotation}\noi}
\nek{\eqo}{\end{quotation}}
\nek{\bpf} {\begin{prF}} 
\nek{\epf} {\qed\end{prF}} 
\nek{\ePf} {\end{prF}} 
\nek{\bpc} {\begin{prC}}
\nek{\epc} {\qeDD{Утверждение}\end{prC}} 
\nek{\qeD} {\qed}
\nek{\qeDD} [1] 
{\hfill\hbox{\qed~({\small #1\/}\hspace{0.1ex})}}
\nek{\epF} [1] {\qeDD{#1}\end{prF}} 

\nek{\bus}{\begin{equation}}   
\nek{\eus}{\end{equation}}
\nek{\bde}{\begin{description}\itsep}
\nek{\ede}{\end{description}}
\nek{\ben}{\begin{enumerate}\itsep}
\nek{\een}{\end{enumerate}}
\nek{\bit}{\begin{itemize}\itsep}
\nek{\eit}{\end{itemize}}
\nek{\bay}{\begin{array}}
\nek{\eay}{\end{array}}
\nek{\bmp}{\begin{minipage}}
\nek{\emp}{\end{minipage}}
\nek{\btb}{\begin{tabular*}}
\nek{\etb}{\end{tabular*}}
\nek{\btab}{\begin{tabular}}
\nek{\etab}{\end{tabular}}
\nek{\fF} {\mathbf F}
\nek{\fG} {\mathbf G}
\nek{\Gd} {\fG_\da}
\nek{\Fs} {\fF_\sg}
\nek{\fS} {{\bf S}}
\nek{\fH} {{\bf H}}
\nek{\fa} {{\bf a}}
\nek{\fb} {{\bf b}}
\nek{\ZFC}{{\bf ZFC}}
\nek{\ZF}{{\bf ZF}}
\nek{\zfcm}{\ZFC^-}
\nek{\ZC} {{\bf ZC}}
\nek{\uqi}{{\bf I}}
\nek{\uqj}{{\phantom{I}\bf I}}
\nek{\uqd}{{\bf II}}

\nek{\gai} [2] {\fH^\cZ_{#1#2}}
\nek{\gad} [2] {\fG^{#1}_{#2}} %{\fG^\cZ_{#1#2}}
\nek{\iai} [1] {\fH_{#1}}
\nek{\iad} [1] {\fG_{#1}}
\nek{\hi} {\fH}
\nek{\hd} {\fG}
\nek{\etc} {{\sl etc}}
\nek{\iesp}{\hspace{0.3ex}}
\nek{\resp}{\hspace{0.25ex}}
\nek{\pw} {\hbox{a.\iesp e.}}
\nek{\ie} {\hbox{\sl i.\iesp e.}}
\nek{\eg} {\hbox{\sl e.\iesp g.}}
\nek{\ea} {\hbox{\sl et.\hspace{0.3ex}al.}}
\nek{\noo}
{\hbox{\sl w.\iesp l.\iesp o.\iesp g.}}
\nek{\pv} [1] {\hbox{почти все{#1}}}
\nek{\vrt} {\hbox{w.\iesp r.\iesp t.}}
\nek{\lsc} {\hbox{l.\iesp s.\iesp c.}}

\nek{\te} {\hbox{т.\resp е.}}
\nek{\td} {\hbox{т.\resp д.}}
\nek{\tp} {\mbox{т.\resp п.\/}}
\nek{\Te} {\hbox{Т.\resp е.}}

\nek{\bore} {\text{\tt BOREQ}}

\nek{\gp} [2] {\text{\tt Gen}^{#1}_{#2}}
\nek{\er} {{ОЭ}}
\nek{\eqr} [1] {отношени{#1} эквивалентности}
\nek{\Eqr} [1] {Отношени{#1} эквивалентности}

\nek{\bm} {{BM}}
\nek{\ddd}[1]{$\mtho\hspace{0.2ex}{#1}$-\hspace{0.0ex}}
\nek{\dw}{\ddd\om}
\nek{\lis}[1] {\mathop{\tt lim\hspace{0.2ex}sup}_{#1}}
\nek{\len}[1] {\mathop{\tt lh}{#1}}
\nek{\Ord}  {{\tt{Ord}}}
\nek{\Exh}  {{\tt{Exh}}}
\nek{\Nul}  {{\tt{Null}}}
\nek{\Mod}  {\mathop{\tt{Mod}}}
\nek{\Aut}  {\mathop{\tt{Aut}}}
\nek{\card} {\mathop{\tt card}}
\nek{\lh}   {\mathop{\tt lh}}
\nek{\pr}   {\mathop{\tt pr}}
\nek{\sr}   {\mathop{\tt sr}}
\nek{\vnr}  {\mathop{\tt rnk}}
\nek{\ran}  {\mathop{\tt ran}}
\nek{\dom}  {\mathop{\tt dom}}
\nek{\fld}  {\mathop{\tt field}}
\nek{\otp}  {\mathop{\tt otp}}
\nek{\Max}  {\mathop{\tt Max}}
\nek{\tsup} {\mathop{\tt sup}}
\nek{\tinf} {\mathop{\tt inf}}
\nek{\tmin} {\mathop{\tt min}}
\nek{\tmax} {\mathop{\tt max}}
\nek{\tlim} {\mathop{\tt lim}}
\nek{\tlis} {\mathop{\tt lim\hspace{0.3ex}sup}}
\nek{\tlii} {\mathop{\tt lim\hspace{0.3ex}inf}}
\nek{\Fin}  {{\tt Fin}} 
\nek{\bFin} {{\bf Fin}} 
\nek{\maxi}[1] {\Max^\xi_{#1}}
\nek{\HC} {{\rm HC}}
\nek{\hc} {\HC}
\nek{\ccc}{{\sc ccc}}
\nek{\al} {\alpha}
\nek{\ba} {\beta}
\nek{\ga} {\gamma}
\nek{\Da} {\Delta}
\nek{\da} {\delta}
\nek{\ka} {\kappa}
\nek{\la} {\lambda}
\nek{\La}{\Lambda}
\nek{\sg} {\sigma}
\nek{\Sg} {\Sigma}
\nek{\vpi}{\varphi}
\nek{\vpy}{\vpi_\iy}
\nek{\ny}{\nu_\iy}
\nek{\vt} {\vartheta}
\nek{\vT} {\Theta}
\nek{\ovt}{{\overline\vt}}
\nek{\ovi}{{\overline\vpi}}
\nek{\ops}{{\overline\psi}}
\nek{\ve} {\varepsilon}
\nek{\om} {\omega}
\nek{\Om} {\Omega}
\nek{\lom}{^{<\om}}
\nek{\za}{\zeta}
\nek{\tpi}{\tau_\vpi}
\nek{\omi} {\om_1}
\nek{\ali} {\aleph_1}
\nek{\omd} {{\om_2}}
\nek{\omm} [1] {\om^{\om^{#1}}}
\nek{\bse} {2\lom} %{\fuk2{<\om}}
\nek{\nse} {\dN\lom}
\nek{\alo} {{\aleph_0}}
\nek{\sd}   {\mathbin{\Da}}
\nek{\bigd} {\mathbin{\hbox{\large$\mtho\Da$}}}
\nek{\sqe} {\sq^\ast}
\nek{\fs}[2]{{\bf\iSg}^{#1}_{#2}}
\nek{\fp}[2]{{\bf\iPi}^{#1}_{#2}}
\nek{\fd}[2]{{\bf\iDa}^{#1}_{#2}}
\nek{\iSg}{{\mathchar"7106}}
\nek{\iPi}{{\mathchar"7105}}
\nek{\iDa}{{\mathchar"7101}}
\newcommand{\is}[2]{\iSg^{#1}_{#2}}

\newcommand{\id}[2]{\iDa^{#1}_{#2}}
\nek{\BBB}{\hspace{0.01ex}}
\nek{\dA}{{\BBB{\mathbb A}\BBB}}
\nek{\dC}{{\BBB{\mathbb C}\BBB}}
\nek{\dF}{{\BBB{\mathbb F}\BBB}}
\nek{\dI}{{\BBB{\mathbb I}\BBB}}
\nek{\dN}{{\BBB{\mathbb N}\BBB}}
\nek{\dP}{{\BBB{\mathbb P}\BBB}}
\nek{\dQ}{{\BBB{\mathbb Q}\BBB}}
\nek{\dqp}{\dQ^+}
\nek{\dR}{{\BBB{\mathbb R}\BBB}}
\nek{\dS}{{\BBB{\mathbb S}\BBB}}
\nek{\dT}{{\BBB{\mathbb T}\BBB}}
\nek{\dV}{{\BBB{\mathbb V}\BBB}}
\nek{\dZ}{{\BBB{\mathbb Z}\BBB}}
\nek{\dX}{{\BBB{\mathbb X}\BBB}}
\nek{\dY}{{\BBB{\mathbb Y}\BBB}}
\nek{\dyn} {\dY^\dN}
\nek{\dvp} {\dV^+}
\nek{\dn}{2^\dN}
\nek{\dntn}{2^{\dN\ti\dN}}
\nek{\dnqn}{{(2^\dN)}{}^\dN}
\nek{\pnqn}{{\pn}{}^\dN}
\nek{\bn}{{\dN\hspace*{-1\mathsurround}}^\dN}
\nek{\rn} {\dR^\dN}
\nek{\rtn}{\rn}
\nek{\ntn}{\bn}
\nek{\ztn}{\dZ^\dN}

\nek{\tm} [1] {2^{\mxi\xi}}
\nek{\nn}{{\dN\ti\dN}}
\nek{\ccs} {}
\nek{\cA}{{\ccs{\skr A}\ccs}}
\nek{\cC}{{\ccs{\skr C}\ccs}}
\nek{\cD}{{\ccs{\skr D}\ccs}}
\nek{\cE}{{\ccs{\skr E}\ccs}}
\nek{\cF}{{\ccs{\skr F}\ccs}}
\nek{\cS}{{\ccs{\skr S}\ccs}}
\nek{\cP}{{\ccs{\skr P}\ccs}}
\nek{\cW}{{\ccs{\skr W}\ccs}}
\nek{\cX}{{\ccs{\skr X}\ccs}}
\nek{\cY}{{\ccs{\skr Y}\ccs}}
\nek{\cI} {{\skr I}} 
\nek{\cJ} {{\skr J}} 
\nek{\cK} {{\skr K}} 
\nek{\cO} {{\skr O}} 
\nek{\cT} {{\skr T}} 
\nek{\cZ} {{\skr Z}}
\nek{\zo} {\cZ_0}
\nek{\zw} {\cZ_{\hbox{\small\rm w}}}
\nek{\xn} {\cX^\dN}
\nek{\an} {\dA^\dN}
\nek{\cd} [1] {\cD_{#1}}
\nek{\pws}  [1] {\cP(#1)}
\nek{\cp}  [2] {\cP_{\hspace*{-0.4ex}\tt cnt}^{#1}(#2)}
\nek{\pwf} [1] {\cP_{\hspace*{-0.4ex}\tt fin}(#1)}
\nek{\pwc} [1] {\cP_{\hspace*{-0.4ex}\tt ctbl}(#1)}
\nek{\pnn}{\pws{\nn}}
\nek{\pn}{\pws{\dN}}
\nek{\ps}{\pws{\dS}}
\nek{\pz}{\pws{\dZ}}
\nek{\mm}{{\BBB{\mathfrak M}\BBB}}
\nek{\mn}{{\BBB{\mathfrak N}\BBB}}
\nek{\gE}{{\BBB{\mathfrak E}\BBB}}
\nek{\ya}  {{\mathfrak a}}
\nek{\shi} {{\mathfrak s}}
\nek{\kon} {\hbox{\mtho\large${\mathfrak c}$}}

\newcommand{\gc}{{\BBB{\mathfrak c}\BBB}}
\nek{\ilo}[1] {{[0\hspace{0.1ex},\hspace{0.1ex}n_{#1})}}
\nek{\ii} [1] {{[n_{#1}\hspace{0.1ex},\hspace{0.1ex}\infty)}}
\nek{\iv} [2] {{(#1\hspace{0.1ex},\hspace{0.1ex}#2)}}
\nek{\ix} [2] {{[#1\hspace{0.1ex},\hspace{0.1ex}#2]}}
\nek{\ir} [2] {{[#1\hspace{0.1ex},\hspace{0.1ex}#2)}}
\nek{\iry}[1] {{[#1\hspace{0.1ex},\hspace{0.1ex}\iy)}}
\nek{\iq} [2] {\ir{\nu_{#1}}{\nu_{#2}}}
\nek{\iqy}[1] {\ir{\nu_{#1}}\iy}
\nek{\iqo}[1] {\ir0{\nu_{#1}}}
\nek{\iqn}[1] {\iq{#1}{#1+1}}
\nek{\opl} {\oplus}
\nek{\ap}  {\cdot}
\nek{\cj}  {\mathbin{\hspace{0.2ex}\&\hspace{0.2ex}}}
\nek{\dm}  {$$}
\nek{\sus} {\mathopen{\exists\hspace{0.35ex}}}
\nek{\kaz} {\mathopen{\forall\hspace{0.35ex}}}
\nek{\imp} {\Longrightarrow} 
\nek{\mpi} {\Longleftarrow} 
\nek{\eqv} {\Longleftrightarrow} 
\nek{\ti}  {\times} 
\nek{\mo}  {\models} 
\nek{\sq}  {\subseteq}
\nek{\qs}  {\supseteq}
\nek{\su}  {\subset}
\nek{\sneq}{\subsetneqq}
\nek{\we}  {{\mathbin{\hspace{0.15ex}^\wedge}}}
\nek{\obr} {^{-1}}
\nek{\dif} {\smallsetminus}
\nek{\res} {\mathbin{\restriction}}
\nek{\lef} {\preccurlyeq}
\nek{\gef} {\succcurlyeq}
\nek{\pu}  {\varnothing}
\nek{\iy}  {\infty}
\nek{\piy} {+\iy}
\nek{\nin} {\not\in}
\nek{\mto} {\mapsto}
\nek{\limp}{\,\imp\,}
\nek{\leqv}{\,\eqv\,}
\nek{\onto}{\stackrel{\text{\rm на}}{\longrightarrow}}
\nek{\cei} [1] {\lceil #1\rceil}
\nek{\ang} [1] {\langle #1\rangle}
\nek{\stk} [2] {\ang{#1\hspace{0.3ex};\hspace{0.1ex}#2}}
\nek{\sis} [2] {\ans{#1}_{#2}}
\nek{\ans} [1] {\{\hspace{0.01ex}#1\hspace{0.01ex}\}}
\nek{\Ans} [1] {\left\{\hspace{0.01ex}#1\hspace{0.01ex}\right\}}

\nek{\zz} {\linebreak[0]} 
\nek{\ens} [2] {\ans{{#1\hspace{0.5ex}{:}}\zz\hspace{0.5ex}#2}}
\nek{\Ens} [2] {\Ans{{#1\hspace{0.5ex}{:}}\zz\hspace{0.5ex}#2}}

\nek{\suh} [1] {[\hspace{0.3pt}#1\hspace{0.3pt}]_{\sq}}
\nek{\itla} {\item\label}

\nek{\aeq} {\mathbin{\|}}
\nek{\laeq}{\,\aeq\,} 

\nek{\tz} {\mathbin{;}}

\nek{\seq}[2] {(#1)_{#2}}

\nek{\kb}[2]{#1^{(#2)}}

\nek{\ssty} {\textstyle}
\nek{\isum} [3] {{{\ssty\ugl\sum_{#1}\,#3}\mid{#2}\ugr}}

\nek{\ifi}  {{\cur{Fin}\hspace*{0.1ex}}}
\nek{\frt}  {{\cur{Fr}\hspace*{0.1ex}}}
\nek{\ibo}  {{\cur{Bou}\hspace*{0.1ex}}}
\nek{\fifi} {\ifi\ti\ifi}
\nek{\fio} {\ifi\ti0}
\nek{\ofi} {0\ti\ifi}

\nek{\xip} {{\xi+1}}
\nek{\etp} {{\eta+1}}

\nek{\vx} [1] {^{(#1)}}

\nek{\dop}[1] {{#1}^{\complement}}

\nek{\skl} {\hbox{\mtho\Large$($}}
\nek{\skp} {\hbox{\mtho\Large$)$}}

\nek{\ugl} {\hbox{\mtho\large$\langle$}}
\nek{\ugr} {\hbox{\mtho\large$\rangle$}}

\nek{\df} [1] {\dop {#1}}
\nek{\pl} [1] {{#1}^+}

\nek{\bbW} {\hbox{\mtho\boldmath $W$}}

\nek{\bo}{{\bf O}}

\nek{\rC}  {\mathbin{\sf C}}
\nek{\rD}  {\mathbin{\sf D}}
\nek{\rE}  {\mathbin{\sf E}}
\nek{\rF}  {\mathbin{\sf F}}
\nek{\rG}  {\mathbin{\sf G}}
\nek{\rM}  {\mathbin{\sf M}}
\nek{\rP}  {\mathbin{\sf P}}
\nek{\rQ}  {\mathopen{{\mathsf Q}\hspace{0.35ex}}}
\nek{\rR}  {\mathbin{\sf R}}
\nek{\rS}  {\mathbin{\sf S}}
\nek{\rT}  {\mathbin{\sf T}}
\nek{\rtt} {\mathbin{\sf t}}
\nek{\rV}  {\mathbin{\tt Vit}}
\nek{\ev}  {\rV}
\nek{\rZ}  {\mathbin{\sf Z}}

\nek{\Zo}  {{\cZ_0}}
\nek{\Eo}  {\rE_{\text{\sf0}}}
\nek{\eow}  {\Eo^{(W)}}
\nek{\eon}  {\Eo^{(\dN)}}
\nek{\eobs} {\Eo^{(\bse)}}
\nek{\Fo}  {\rF_0}
\nek{\Ec}  {\mathbin{\overline{\rE}}}
\nek{\Eco} {\mathbin{\overline{\Eo}}}
\nek{\rzo}  {\rZ_{\text{\sf0}}}
\nek{\neo} {\mathbin{\not{{\hspace{-0.4ex}\rE}}_0}}
\nek{\nfo} {\mathbin{\not{{\hspace{-0.4ex}\rF}}_0}}
\nek{\nrE} [1] {\mathbin{\not{{\hspace{-0.4ex}\rE}}_{#1}}}
\nek{\rtd} {\rT_2}
\nek{\rdi} {\mathbin{{\sf D}_I}}
\nek{\rda} {\mathbin{{\sf D}_A}}

\nek{\dno} {\not\hspace*{0.05ex}\not}
\nek{\scp} {{\text{\bf+}}}
\nek{\epo} [1] {\mathbin{{#1}{\vphantom{\pn}}^\scp}}
\nek{\ei}   {\epo{\rE}}
\nek{\rfi}  {\epo{\rF}}
\nek{\nE}  {\mathbin{{\dno\hspace{-0.35ex}\sf E}}}
\nek{\nF}  {\mathbin{{\dno\hspace{-0.25ex}\sf F}}}
\nek{\reo} {\rE_{\hspace{-1.0pt}\Zo}}
\nek{\rez} {\rE_{\hspace{-1.0pt}\cZ}}
\nek{\nrz} {\mathbin
{\not{{\hspace{-0.4ex}\rE}}_{\hspace{-1.0pt}\cZ}}}
\nek{\nre} {\mathbin
{\not{{\hspace{-0.4ex}\rE}}_{\hspace{-1.0pt}\Zo}}}
\nek{\reff}{\rE_{\fifi}}
\nek{\dzo} {\ddd{\Zo}}
\nek{\dde} {\ddd{\rE}}
\nek{\ddec}{\ddd{\Ec}}
\nek{\dep} {\ddd{\rE'}}
\nek{\ddf} {\ddd{\rF}}
\nek{\dfp} {\ddd{\rF'}}
\nek{\ddv} {\ddd{\vt}}
\nek{\ddz} {\ddd\cZ} 
\nek{\der} {\er-}
\nek{\dee} {\ddd{\rE,\rE'}}
\nek{\Def} {\ddd{\rE,\rF}}

\nek{\fdo}{\hbox{\raisebox{0.2ex}{\mtho\tiny$\bullet$}}}
\nek{\fdt}{\hbox{\raisebox{-0.25ex}{\LARGE\bf.}}}
\nek{\bdot}[1] {\raisebox{-0.07ex}{\mtho$\stackrel{\fdt}{#1}$}}
\nek{\dvpi} {{\bdot\vpi}} 
\nek{\doal} {{\bdot\al}} 
\nek{\doga} {{\bdot\ga}} 
\nek{\doa} {{\bdot a}} 
\nek{\don} {{\bdot n}} 
\nek{\dotm} {{\bdot m}} 
\nek{\dok} {{\bdot k}} 
\nek{\doll} {{\bdot\ell}} 
\nek{\dox} {{\bdot x}}
\nek{\doz} {{\bdot z}}
\nek{\doy} {\raisebox{-0.28ex}{\mtho$\stackrel{\fdt}y$}}
\nek{\dtp} {\raisebox{-0.28ex}{\mtho$\stackrel{\fdt}p$}}
\nek{\dog} {\text{\bf g}}
\nek{\doxl}{\dox_{\tt left}}
\nek{\doxr}{\dox_{\tt right}}
\nek{\rkb} [1] {|#1|_{\rm CB}}
\nek{\rkt} [2] {{^{\om}\hspace*{-1.3pt}|#1|_{#2}}}
\nek{\rko} [2] {{|#1|_{#2}}}
\nek{\rkT} [1] {{^{\om}\hspace*{-1.3pt}|#1|}}
\nek{\rkO} [1] {{|#1|}}
\nek{\kt}  [1] {{^{#1}\hspace*{-0.8pt}T}}
\nek{\bsf} [1] {I_{#1}}
\nek{\fri} [1] {\mathbin{\displaystyle\rE^{\tt fr}_{#1}}}
\nek{\fre} [1] {\frt_{{#1}}}
\nek{\fss} [3] {{\sum_{#3}{#2}\,/\,{#1}}}
\nek{\fps} [3] {{\prod_{#3}{#2}\,/\,{#1}}}
\nek{\fpt} [2] {{\prod{\sis{#2}{}}\,/\,{#1}}}
\nek{\fpd} [2] {{#1}\otimes{#2}}
\nek{\fsm} [2] {{#1}\oplus{#2}}

\nek{\sto} {[s_0,t_0]}
\nek{\ostp}{[s',t']}
\nek{\ost} {[s,t]}
\nek{\osu} {[s,u]} 
\nek{\otu} {[t,u]}
\nek{\ouv} {[u,v]}

\nek{\zs} {{\tilde s}}
\nek{\zt} {{\tilde t}}
\nek{\zu} {{\tilde u}}
\nek{\zv} {{\tilde v}}
\nek{\zn} {{\tilde n}}
\nek{\zm} {{\tilde m}}
\nek{\zx} {{\tilde x}}
\nek{\zU} {{\widetilde U}}
\nek{\zI} {{\widetilde \cI}}

\nek{\ff}[2] {F_{#1}^{#2}}

\nek{\spa} {\dS}
\nek{\spn} {\spa^\dN}
\nek{\spp} {\dY}

\nek{\poq} {\underline}
\nek{\nad} {\overline}
\nek{\nadd}{\widehat}

\nek{\sm} {\hbox{\mtho{\large$\Sg$}}}
\nek{\smp} [2] {\sm_{#1}^{#2-1}}
\nek{\smy} [1] {\sm_{#1}^{\iy}}

\nek{\sui} [1] {\cS_{\ans{#1}}}
\nek{\sun} {{\sui{\frac1{n+1}}}}
\nek{\srn} {{\sui{r_n}}}
\nek{\ern} {\rS_{\ans{r_n}}}
\nek{\erpn} {\rS_{\ans{r'_n}}}
\nek{\eun} {\rS_{\ans{\frac1{n+1}}}}
\nek{\suo} {{\sui0}}
\nek{\eso} {\rE_{\hspace{-1.0pt}\suo}}
\nek{\nso} {\mathbin
{\not{{\hspace{-0.4ex}\rE}}_{\hspace{-1.0pt}\suo}}}

\nek{\gal} [3] {{\tt Gal}^{#1}_{#2}(#3)} 

\nek{\nbd} [1] {{\skr O}_1(#1)}

\nek{\aH} {H^\ast}
\nek{\aB} {B^\ast}

\nek{\ovl} [1] {\overline{#1}} 
\nek{\ovg} [1] {\ovl{g_{#1}}}
\nek{\ovp} [1] {\ovl{\ga_{#1}}}

\nek{\plo} {+1} %{\!+\!1}

\nek{\yk} [1] {k_{#1}}

\nek{\mtho}{\mathsurround=0mm}
\nek{\msur}{\hspace{-1\mathsurround}}
\nek{\psur}{\hspace{0.3\mathsurround}}
\nek{\dsur}{\hspace{-0.3\mathsurround}}
\nek{\hsur}{\hspace{-0.5\mathsurround}}
\nek{\noi}{\noindent}
\nek{\vom}{\vspace{1mm}}
\nek{\vtm}{\vspace{2mm}}

\nek{\uv}{{\bf V}}

\nek{\wA} {{\widehat A}}
\nek{\ha} {{\hat a}}
\nek{\hb} {{\hat b}}
\nek{\he} {{\hat\ve}}
\nek{\hT} {{\hat t}}
\nek{\hl} {{\hat l}}
\nek{\hs} {{\hat s}}
\nek{\hm} {{\widehat m}}
\nek{\hx} {{\hat x}}
\nek{\hy} {{\hat y}}
\nek{\hf} {{\widehat f}}
\nek{\hn} {{\widehat n}}
\nek{\hk} {{\widehat k}}
\nek{\hvt}{{\widehat\vt}}

\nek{\fo} {{\mathbf 0}}
\nek{\fr} {{\mathbf 1}}

\nek{\vnu} {{\vec \nu}} 
\nek{\dvn} {\ddd\vnu} 
\nek{\wnu} {\cW_{\vnu}} 
\nek{\ewn} {\rE_{\vnu}}
\nek{\ret} {\rE_{T}}

\nek{\ci} [1] {I_{#1}}

\nek{\vex}{{\vec x}}
\nek{\vey}{{\vec y}}

\nek{\dns} [1] {d_{\vnu}^{#1}}
\nek{\den} {d_{\vnu}}
\nek{\din} {d_\vnu}

\nek{\poo}{=_{\tt df}}

\nek{\dpi}{d_\vpi}

\nek{\nol}[1] {{\bf 0}_{#1}}
\nek{\edi}[1] {{\bf 1}_{#1}}

\nek{\nrn} {_{\ans{r_n}}}
\nek{\vpr} {\vpi\nrn}
\nek{\dpr} {d\nrn}
\nek{\drn} {\ddd{\ans{r_n}}}

\nek{\okr} [2] {{\skr O}_{#1}(#2)}

\nek{\nid}{\gE}

\nek{\zid} [2] {\cI_{#1/#2}}
\nek{\zidef} {\zid\rE\rF}

\nek{\zfo} [2] {\dP_{#1/#2}}
\nek{\zfoef} {\zfo\rE\rF}

\nek{\peo} {\dP_{\Eo}}
\nek{\pep} {\dP'_{\Eo}}
\nek{\ieo} {\cI_{\Eo}}

\SetMathAlphabet{\cur}{bold}{U}{eur}{b}{n}
\renek{\Gd} {\fG_\fda}
\renek{\Fs} {{\fF\hspace*{-0.2ex}}_\fsg}

\nek{\ubf}{\fontseries{b}\selectfont}

\mathchardef\alphA ="710B
\mathchardef\betA ="710C
\mathchardef\gammA ="710D
\mathchardef\deltA ="710E
\mathchardef\vartA ="7123
\mathchardef\kpA   ="7114
\mathchardef\lA    ="7115
\mathchardef\mU    ="7116
\mathchardef\nU    ="7117
\mathchardef\rhO   ="711A
\mathchardef\sigmA ="711B

\nek{\fal} {{\cur{\alphA}}}
\nek{\fba} {{\cur{\betA}}}
\nek{\fpi} {{\cur{\pI}}}
\nek{\fsg} {{\cur{\sigmA}}}
\nek{\fda} {{\cur{\deltA}}}
\nek{\fla} {{\cur{\lA}}}

\nek{\dds}{\ddd\fsg}

\nek{\nmp} {\Longleftarrow}

\nek{\tsc}[1]{\hbox{\footnotesize\sc{#1}}}

\nek{\ddi} {\ddd{\cI}}
\nek{\ddj} {\ddd{\cJ}}

\nek{\ren} {\le_{\tsc{c}}}
\nek{\rens} {<_{\tsc{c}}}
\nek{\eqN} {\sim_{\tsc{c}}}

\nek{\eui}[1] {{\text{\it{EU}}}_{\ans{#1}}}

\nek{\sqa} {\sq^\ast}
\nek{\sua} {\su^\ast}

\nek{\ppp} [1] {\hbox{пп.\hspace{0.3ex}\ref{#1}}}
\nek{\ppf} [1] {\hbox{п.\hspace{0.3ex}\ref{#1}}}
\nek{\prf} [1] {\hbox{\S\hspace{0.3ex}\ref{#1}}}
\nek{\prff}[1] {\S\S\hspace{0.3ex}\ref{#1}}
\nek{\pff} [1] {\S\:{#1}}

\nek{\mrn} {\mu\nrn}

\renek{\dop} [1] {\complement #1}
\renek{\df}  [1] {{#1}^\complement}
\nek{\doP}  [1] {{#1}^\complement}

\nek{\ima} [2] {{#1}[#2]}
\nek{\aprb}{\approx_{\tsc b}}
\nek{\ismb}{\cong_{\tsc b}}

\nek{\reas} {<_{\tsc a}}
\nek{\rea} {\leqslant_{\tsc a}}
\nek{\eqa} {\sim_{\tsc a}}

\nek{\geb} {\geqslant_{\tsc b}}
\nek{\rebs} {<_{\tsc b}}
\nek{\reb} {\leqslant_{\tsc b}}
\nek{\eqb} {\sim_{\tsc b}}
\nek{\izb} {\cong_{\tsc b}}

\nek{\reBs}{<_{\tsc{bm}}}
\nek{\reB} {\leqslant_{\tsc{bm}}}
\nek{\eqB} {\sim_{\tsc{bm}}}

\nek{\rds} {<_{\tsc{b},\Da}}
\nek{\rd} {\leqslant_{\tsc{b},\Da}}
\nek{\eqd} {\approx_{\tsc{b},\Da}}

\nek{\rbas} {<_{\tsc{b,ba}}}
\nek{\rba} {\leqslant_{\tsc{b,ba}}}
\nek{\eqba} {\approx_{\tsc{b,ba}}}

\nek{\rbasp} {<_{\tsc{b,ba}}^+}
\nek{\rbap} {\leqslant_{\tsc{b,ba}}^+}
\nek{\eqbab} {\approx_{\tsc{b,ba}}^+}

\nek{\nab} [1] {\nabla(#1)}
\nek{\orb} {\leqslant_{\tsc{rb}}}
\nek{\srb} {<_{\tsc{rb}}}
\nek{\erb} {\sim_{\tsc{rb}}}

\nek{\orbpp} {\leqslant_{\tsc{rb}}^{++}}
\nek{\srbpp} {<_{\tsc{rb}}^{++}}
\nek{\erbpp} {\sim_{\tsc{rb}}^{++}}

\nek{\orbp} {\leqslant_{\tsc{rb}}^{+}}
\nek{\srbp} {<_{\tsc{rb}}^{+}}
\nek{\erbp} {\sim_{\tsc{rb}}^{+}}

\nek{\ork} {\leqslant_{\tsc{rk}}}
\nek{\srk} {<_{\tsc{rk}}}
\nek{\erk} {\sim_{\tsc{rk}}}

\nek{\obe} {\leqslant_{\tsc{be}}}
\nek{\sbe} {<_{\tsc{be}}}
\nek{\ebe} {\sim_{\tsc{be}}}

\nek{\odl} {\leqslant_{\sd}}
\nek{\sdl} {<_{\sd}}
\nek{\edl} {\sim_{\sd}}

\nek{\rei} {\rE_{\cI}}
\nek{\rej} {\rE_{\cJ}}

\nek{\eeb} {\sim_{\tsc{b}}}

\nek{\obep} {\leqslant_{\tsc{be}}^+}
\nek{\sbep} {<_{\tsc{be}}^+}
\nek{\ebep} {\sim_{\tsc{be}}^+}

\nek{\seb} {<_{\tsc{b}}}

\nek{\supp} {\mathop{\tt supp}}

\nek{\atp} {\mathop{\tt at}^+}
\nek{\atm} {\mathop{\tt at}^-}

\nek{\hv} [2] {||#1||_{#2}}

\nek{\vmu} {{\vec \mu}}

\nek{\bel} [1] {\mathbin{\text{\boldmath\mtho$\ell$}}^{#1}}
\nek{\beL} [1] {\mathbin{\text{\bfsf L}}^{#1}}
\nek{\bem}     {\mathbin{\text{\bfsf m}}}
\nek{\fCo}{\mathbin{\rC_{\text{\sf0}}}}
\nek{\fco}{\mathbin{\text{\bfsf c}_{\text{\sf0}}}}
\nek{\fc} {\mathbin{\text{\bfsf c}}}
\nek{\fvt}{\mathbin{\text{\boldmath\mtho$\vt$}}}

\nek{\beli} {\bel\iy}

\nek{\oin} {{[0,1]^\dN}}

\nek{\iz} [2] {\cI_{#1}^{#2}}
\nek{\jz} [1] {\iz{}{}(#1)}
\nek{\iw} [2] {\cW_{#1}^{#2}}
\nek{\jw} [1] {\iw{}{}(#1)}
\nek{\ib} [2] {\cB^{#1}_{#2}}
\nek{\jb} [1] {\ib{}{#1}}

\nek{\ovu} {{\overline u}}
\nek{\ovv} {{\overline v}}

\nek{\nr}[2]{\nor{#1}_{#2}} %{||#1||_{#2}}

\nek{\TS}{\textstyle}
\nek{\DS}{\displaystyle}

\nek{\Ba}{B^\ast}

\nek{\resi} [1] {\mathop{\restriction_{#1}}}

\nek{\gP}{{\BBB{\mathfrak P}\BBB}}
\nek{\gF}{{\BBB{\mathfrak F}\BBB}}
\nek{\gJ}{{\BBB{\mathfrak J}\BBB}}
\nek{\dG}{{\BBB{\mathbb G}\BBB}}
\nek{\dH}{{\BBB{\mathbb H}\BBB}}

\nek{\ac} {\cdot} %action

\nek{\curle}{\preccurlyeq}
\nek{\cle}{\curle}
\nek{\cge}{\succcurlyeq}
\nek{\cl} {\prec}
\nek{\resic} [1] {\resi{\cle #1}} 
\nek{\rec} {\resic} 
\nek{\recb} [1] {\rec{(#1)}} 
\nek{\rsic} [1] {\resi{\cl #1}} 
\nek{\rc} {\rsic} 
\nek{\rcb} [1] {\rc{(#1)}} 

\nek{\kai} {\forall^\iy\hspace{0.1ex}}
\nek{\exi} {\exists^\iy\hspace{0.1ex}}

\nek{\ovw}{\nad w}

\nek{\il}[2] {\ir{n_{#1}}{n_{#2}}}
\nek{\ia}[2] {\ir{a_{#1}}{a_{#2}}}
\renek{\ij}[2] {\ir{j_{#1}}{j_{#2}}}
\nek{\cB}{{\BBB{\skr B}\BBB}}
\nek{\cG}{{\BBB{\skr G}\BBB}}
\nek{\cL}{{\BBB{\skr L}\BBB}}
\nek{\cN}{{\BBB{\skr N}\BBB}}
\nek{\cU}{{\BBB{\skr U}\BBB}}

\nek{\pnd}{\pn^\dN}
\nek{\dnd}{(\dn){}^\dN}

\nek{\anp} [2] {\ang{#1}^{#2}}

\nek{\ta}{\tau}

\nek{\lev} {\mathop{\tt{lev}}}
\nek{\glu} {\mathop{\tt{dep}}}
\nek{\dia} {\mathop{\tt{diam\hspace{0.15ex}}}}

\nek{\Ei}{\rE_{\text{\sf 1}}}
\nek{\Ed}{\rE_{\text{\sf 2}}}
\nek{\Et}{\rE_{\text{\sf 3}}}
\nek{\Eh}{\rE_{\text{\sf t}}}
\nek{\Ey}{\rE_\iy}

\nek{\npi}{\nu_\vpi}
\nek{\nsi}{\nu_\psi}

\nek{\Ii}{\cI_1}
\nek{\Id}{\cI_2}
\nek{\It}{\cI_3}

\nek{\emb}  {\sqsubseteq_{\tsc{b}}}
\nek{\emn}  {\sqsubseteq_{\tsc{c}}}
\nek{\embi} {\sqsubseteq_{\tsc{b}}^{\rm i}}
\nek{\emni} {\sqsubseteq_{\tsc{c}}^{\rm i}}

\nek{\sio}{\cS_0}
\nek{\esn}{\rE_{\ans{1/n}}}
\nek{\nsn}{\mathbin{\not\hspace*{-0.3ex}\esn}}

\nek{\req}[2]{{\DS|_{#1}^{#2}}}
\nek{\rlq}[1]{\resi{\ge#1}}
\nek{\rmq}[1]{\resi{<#1}}

\nek{\Dij} {\ddd{\cI,\cJ}}
\nek{\deo} {\ddd{\Eo}}

\nek{\dc}[2] {{\bf W}^{#1}_{#2}}

\nek{\Za}{{\nad Q}}
\nek{\Ya}{{\widetilde Y}}

\nek{\Ga}{\Gamma}

\nek{\inva}{{\tt{inv}}}

\nek{\tP}{{P^\ast}}
\nek{\app} {{\hspace{0.2ex}{\cdot}\hspace{0.2ex}}}

\nek{\lland}{\,\land\,}

\nek{\seqv} {\hspace{0.3ex}\Leftrightarrow\hspace{0.3ex}}
\nek{\simp} {\hspace{0.3ex}\Rightarrow\hspace{0.3ex}}

\nek{\eqn}{\equiv_n}

\nek{\pit}{\tilde\pi}
\nek{\tve}{\tilde\ve}

\nek{\dsu}[2] {{#1\oplus#2}}

\nek{\ef} {\ddd{\rE,\rF}}
\nek{\efp}{\ddd{\rE,{\rF'}}}

\nek{\isi} {\cong}
\nek{\isa} {\mathrel{\hspace{0.2ex}\isi^\ast}}

\nek{\ske} [3] {\equiv_{#1#2}^{#3}}
\nek{\skab}[1] {\ske AB{#1}}
\nek{\spab}[1] {\ske{A'}{B'}{#1}}

\nek{\Indent}{\hspace*{3ex}}
\nek{\fsur}{\hspace{0.5\mathsurround}}

\nek{\rdm}  {\rD_{\tt max}}

\nek{\co} {\ddd{\fco}}
\nek{\lv} {{\normalfont\scshape lv}-}

\nek{\susi} {\exists^{\iy}\hspace*{0.2ex}}
\nek{\kazi} {\forall^{\iy}\hspace*{0.2ex}}

\nek{\igi}{$\hbox{\mtho\boldmath$\fal$}$}
\nek{\igd}{$\hbox{\mtho\boldmath$\fba$}$}

\nek{\cg}[1] {{\tt Choq}(#1)}
\nek{\cgs}[1]{{\tt Choq}^{\rm s}(#1)}

\nek{\dnn}{\dN^\dN}
\nek{\dqq}{\dQ^\dQ}
\nek{\nnn}{(\dnn){\vphantom{\dN}}^\dN}

\nek{\isg} {{S_\iy}}
\nek{\Sy}  {\isg}

\nek{\uset}{{Universal sets}}

\nek{\ler}[2] {\mathbin{\sim^{#1}_{#2}}}
\nek{\aer}[2] {\mathbin{\rE_{#1}^{#2}}}
\nek{\ergx}{\aer\dG\dX}
\nek{\egx} {\ergx}
\nek{\lo} [3] {\cO(#3,#1,#2)}
\nek{\sym}{\ler} % [2] {\sim_{#1}^{#2}}
\nek{\ong} {1_\dG}
\nek{\rr} [2] {\rR^{#1}_{#2}}

\nek{\rav}[1] {{\mathsf\Da}_{#1}}

\nek{\toq}{\circle*{0.5}}
\nek{\tob}{\circle*{1.0}}

\nek{\stob}{\circle{3.0}}
\nek{\ktob}{\kras\circle{3.0}}

\nek{\cob}{\circle{1.5}}
\nek{\mtir}{\line(-1,0){2}}

\nek{\bon} [2] {[#1]}

\nek{\dnnn} {\dnnp\dN}
\nek{\dnnp} [1] {(\dnn){\vphantom{\dN}}^{#1}}

\nek{\prift}{\sf}
\nek{\Penu} [1] {{\prift\ddd{\id11}кодировк{#1}}}
\nek{\Refl} [1] {{\prift Отражени{#1}}}
\nek{\Uset} [1] {{\prift Теорем{#1} Универсальных Множеств}}
\nek{\Kres} [1] {{\prift Выбор{#1} по Крайзелю}}
\nek{\Cenu} [2] {{\prift Счетно-формн{#1} Перечислени{#2}}}
\nek{\Cuni} [2] {{\prift Счетно-формн{#1} Униформизаци{#2}}}
\nek{\Kuni} [2] {{\prift \ddd\sg компактн{#1} Униформизаци{#2}}}
\nek{\Cpro} [2] {{\prift Счетно-формн{#1} Проекци{#2}}}
\nek{\Kpro} [2] {{\prift \ddd\sg компактн{#1} Проекци{#2}}}
\nek{\Sepa} [1] {{\prift Отделимост{#1}}}
\nek{\Redu} [1] {{\prift Редукци{#1}}}
\nek{\Prod} [2] {{\prift Теорем{#1} продолжени{#2}}}

\nek{\dpf} {\mathord{{\dP^2\hspace*{-0.3ex}}\res\rF}}
\nek{\dpe} {\mathord{{\dP^2\hspace*{-0.3ex}}\res\rE}}

\nek{\ek}[2] {[#1]_{{#2}}}

\nek{\eke}[1] {\ek{#1}{\rE}}
\nek{\ekeo}[1] {\ek{#1}{\Eo}}
\nek{\ekec}[1] {\ek{#1}{\Eco}}
\nek{\ekco}[1] {\ek{#1}{\Ec}}
\nek{\ekfo}[1] {\ek{#1}{\Fo}}
\nek{\ekf}[1] {\ek{#1}{\rF}}
\nek{\ekg}[1] {\ek{#1}{G}}

\nek{\ur}{_{\tt right}}
\nek{\ul}{_{\tt left}}

\nek{\cont}{\hbox{\mtho\large$\gc$}}

\nek{\mem} {\ddd\in}

\nek{\bL}{{\bf L}} 
\nek{\cli} {{\sc cli}} 

\nek{\di} [1] {{#1}^\#}
\nek{\drE} {\mathbin{\di\rE}}
\nek{\drF} {\mathbin{\di\rF}}
 
\nek{\bk} [1] {{\cur B}_{#1}}

\nek{\wtau}{{\widehat\tau}}

\nek{\gra}[1]{\ddd{#1}``grainy''}
\nek{\grap}{``grainy''}

\nek{\xE}[2] {\mathbin{\rR^{#2}_{\ge #1}}} 

\nek{\moq} [1] {\Mod_{#1}}
\nek{\mox} {\moq} %[1] {{\mathbb X}_{#1}}
\nek{\loa} [1] {j_{#1}}
\nek{\ism} [1] {\cong_{#1}}
\nek{\izm} [2] {\cong_{#1}^{#2}}
\nek{\aut} [1] {\Aut_{#1}}

\nek{\hfn} {{\rm HF}(\dA)}
\nek{\tce} [1] {{\rm TC}_\ve(#1)}
\nek{\ihf} {\simeq_{\hfn}}

\nek{\symr}{\equiv}
\nek{\rrt} [3] {\symr_{#2#3}^{#1}}
\nek{\rrq} [5] {{#4}\symr_{#2#3}^{#1}{#5}}
\nek{\nrq} [5] {{#4}\not\symr_{#2#3}^{#1}{#5}}
\nek{\rrQ} [5] {{#4}\symr_{#2\,,\,#3}^{#1}{#5}}
\nek{\rro} [5] {\ang{#2,#4}\symr^{#1}\ang{#3,#5}}
\nek{\rrO} [1] {\symr^{#1}}

\nek{\lww} {\cL_{\omi\om}}

\nek{\forc}
{\mathrel{{|\hspace{-0.2ex}|\hspace{-1.1ex}-\hspace{-1.7ex}-}}}

\newlength{\dxii}
\nek{\fC} {{\bf C}}
\nek{\pg} {\fC_\dG}
\nek{\px} {\fC_\dX}
\nek{\pxx} {\fC_{\dX\ti\dX}}
\nek{\pgx} {\fC_{\dG\ti\dX}}

\nek{\nasl} [1] {локальн{#1}}

\nek{\gen} {ген.\ }
\nek{\gene} {ген.}
\nek{\hg}  {лок.\hspace{0.4ex}ген.\ }
\nek{\hgp} {лок.\hspace{0.4ex}ген.}
\nek{\incl} [1] {\mathop{\text{\sc Int}}\overline{#1}}  %{{\tt IntCl}}

\nek{\PP} {pinned}
\nek{\PPP}{Pinned}

\nek{\bap}{{\bar p}}
\renek{\wtau} {{\widehat p}}

\nek{\zO} {\yo}
\nek{\zi} {\yi}
\nek{\zd} {\yd}
\nek{\zT} {\yt}
\nek{\yz} {\linebreak[0]}
\nek{\yo} {,\linebreak[0]}
\nek{\yi} {\hspace{0.2ex},\linebreak[0]\hspace{0.2ex}}
\nek{\yd} {\hspace{0.2ex},\linebreak[0]\:}
\nek{\yt} {\hspace*{0.2ex},\linebreak[0]\;}

\nek{\prit} [1] {[{{\rm #1}}]}
\nek{\nor} [1] {\|#1\|}

\nek{\fap} {f.\hspace{0.1ex}a.\hspace{0.1ex}p.\hspace{0.1ex}m.}

\nek{\bpr} [1] {\bpf[{{\sl#1}\/}]}

\nek{\fx} {{\bf x}}

\nek{\rzd} {^{\text{\tt red}}}

\nek{\srez} [2] {{\bf S}_{#1}(#2)}

\nek{\qc}  [1] {\resi{> #1}} 
\nek{\qec} [1] {\resi{\ge #1}} 
\nek{\rme} [1] {\resi{\le#1}}

\nek{\alex} {<_{\text{\tt alex}}}
\nek{\lex} {<_{\text{\tt lex}}}
\nek{\lexe} {\le_{\text{\tt lex}}}
\nek{\act} {<_{\text{\tt act}}}

\nek{\bbo} {\mathbf O}

\nek{\fz} {{\mathbf z}}

\nek{\dln} {2\lom}

\nek{\fK} {{\bf K}}

\nek{\Ks} {\fK_\fsg}

\nek{\dva}{{\ans{0,1}}}

\nek{\snos} [1] {\,\footnote{\hspace*{2pt}#1}}
\nek{\snom}     {\,\footnotemark}
\nek{\snot} [1] {\,\footnotetext{\hspace*{2pt}#1}}

\nek{\rH}  {\mathbin{\sf H}}

\nek{\fras}[2] {\text{\footnotesize$\DS\frac{#1}{#2}$}}
\nek{\fral}[2] {\text{\large$\frac{#1}{#2}$}}

\nek{\renu}{\tenu{{\rm(\roman{enumi})}}}
\nek{\Renu}{\tenu{{\rm(\Roman{enumi})}}}
\nek{\aenu}{\tenu{{\rm(\arabic{enumi})}}}
\nek{\aenus}{\tenu{{\mtho$(\text{\arabic{enumi}}')$}}}

\nek{\ergg}{\aer\dG\dG}

\nek{\rit} [1] {{\it#1\/}}

\nek{\lap} [1] {\text{<<}#1\text{>>}}

\nek{\fm} [2] {{#1}/{#2}}

\nek{\fuk} [2] {\mathord{^{#2}\hspace*{-0.3ex}{#1}}}

\nek{\rbt} [3] {[#1]^{#2}_{#3}}

\nek{\har} [1] {\chi_{#1}}

\nek{\lam} [1]
{\label{#1}\hspace*{-3pt}\imar{#1}%
}%
\nek{\las} [1]
{\label{#1}\imae{#1}}%
\newcounter{pa}
\newcounter{pb}
\nek{\pai}{%\phantom{1}
\setcounter{pa}{\value{page}}}

\nek{\nat} [1] {\#(#1)}

\nek{\pgcr} {\addtocounter{page}1}

\vyk{
\nek{\pgcr}{\setcounter{pb}{\value{page}}%
\ifodd\value{pa}
\ifodd\value{pb}
\addtocounter{page}1
\fi
\else
\ifodd\value{pb}\else\addtocounter{page}1
\fi
\fi
}
}

\nek{\mgcr} {\addtocounter{page}{-1}}

\nek{\imas}[1]{\marginpar
[\mtho\flushright\footnotesize%
{\vspace{-3ex}$\mtho\longrightarrow$}\\ 
{#1}]
{\mtho\flushleft\footnotesize% 
{\vspace{-3ex}$\mtho\longleftarrow$}\\ 
{#1}}}

\nek{\imad}[1]{\marginpar
[\mtho\flushright\footnotesize%
{\vspace{3ex}$\mtho\longrightarrow$}\\ 
{#1}]
{\mtho\flushleft\footnotesize% 
{\vspace{3ex}$\mtho\longleftarrow$}\\ 
{#1}}}

\nek{\urav}  [1] {\bus\itsep#1\eus} 
\nek{\uravm} [2]
{\imad{#1}\bus{\label{#1}}\itsep#2\eus} 

\nek{\busm} [1] {%
\imad{#1} 
\bus\label{#1}%
}

\nek{\vykl}{\vyk}                                   
\nek{\onon} {\/\text{\rm1\,--\,1}}

\nek{\mast} {{\ensuremath{(\ast)}}}
\nek{\mdag} {{\ensuremath{(\ast)}}}

\nek{\arr} [1]
{\overset{#1}{\longrightarrow}}
\nek{\dnz} {(\dn)^\dZ}

\nek{\inw} [3] {\ddd{(#1\to#2)}инвариантн{#3}}

3

\nek{\uoa} {\mathbin{{\text{\bfsf u}}^\ast_0}}
\nek{\uo}  {\mathbin{{\text{\bfsf u}}_0}}

\nek{\ld} [1] {\mathbin{{\jd}(#1)}}
\nek{\xd} [2] {{\mathbin{{\jd}(#1\hspace*{0.1ex};\hspace*{0.1ex}#2)}}}
\nek{\qd} [3] {{\mathbin{{\jd}
(\ang{#1\hspace*{0.1ex};\hspace*{0.1ex}#2}_{#3\in\dN})}}}
\nek{\qqd} [3] {{\mathbin{{\jd}
(\ang{#1\hspace*{0.1ex};\hspace*{0.1ex}#2}_{#3})}}}

\nek{\jd}  {\mathbin{\rD}}
\nek{\ntd} {\mathbin{{\not\hspace{-0.0ex}\jd}}}

\nek{\nrf}{\prf}
\nek{\gla} [2] {глав{#1}~\ref{#2}}

\nek{\zad} [1] {{\ubf упражнени{#1}}}
\nek{\Zad} [1] {{\ubf Упражнени{#1}}}

\nek{\bom} [2] {{бо\-ре\-лев\-ск{#1} мощ\-ност{#2}}}

\nek{\kom} [1] {{\mathbf K}(#1)}

\nek{\ibn}[1] {\bon{#1}{\bn}}
\nek{\ibd}[1] {\bon{#1}{\dn}}

\nek{\ke}[1]{[#1]} 
   
\nek{\gopa}{\bf}  %  {\GOT}
\nek{\gW}{{\gopa W}}
\nek{\gD}{{\gopa D}}

\nek{\api}{\addtocounter{page}{1}}
\nek{\atc}{\addtocounter{enumi}{1}}
\nek{\cn} [2] {\dN^{#1}\ti {(\bn)}{\vphantom|}^{#2}}
\nek{\xm} [1] {\vhm\dm#1\vhm\dm}
\nek{\vhm}{\vspace{0.5mm}}

\nek{\ns} [1] {{\text{\mtho\boldmath${\cur s}\hspace{0.0ex}$}}_{#1}}

\nek{\bbla}{\begin{bla}}
\nek{\ebla}{\qed\end{bla}}

\nek{\bqu}{\begin{quotation}\noi}
\nek{\equ}{\end{quotation}}

\nek{\Bor} [1] {{\bf{Bor}}(#1)}

\nek{\btB}{\begin{tabular}}
\nek{\etB}{\end{tabular}}

\nek{\mthf}{\mathsurround=0.25ex}

\nek{\qq} [1] {{\sc#1}}
\nek{\na}{\onto}
\nek{\imo} [2] {#1\obr[#2]}

\nek{\iGa}{{\mathchar"7100}}
\nek{\ig}[2]{\iGa^{#1}_{#2}}
\nek{\fg}[2]{{\bf\iGa}^{#1}_{#2}}
\nek{\fGa} {{\mathbf\iGa}}

\nek{\paf} {Параграф}

\nek{\fO} {\text{\ubf 0}}

\nek{\graf} [1] {\Ga_{#1}}
\nek{\grag} [1] {\Ga_{\raz #1}}

\nek{\rk}[2] {|{#2}|_{#1}}
\nek{\rak} [1] {|{#1}|}
\nek{\rai} [1] {|{#1}|_{\text{\tt wf}}}

\nek{\tre} {\text{\ubf Tr}}
\nek{\wft} {\text{\ubf WFT}}
\nek{\ift} {\text{\ubf IFT}}

\nek{\wo} {\text{\ubf WO}}
\nek{\io} {\text{\ubf IO}}

\nek{\hz} {\hat z}

\nek{\pA} {\dA}
\nek{\pX} {\dX}
\nek{\pY} {\dY}
\nek{\pS} {\dS}

\nek{\mov}{\vspace{-1mm}}

\nek{\wR} {\widehat R}
\nek{\wS} {\widehat S}
\nek{\wvt} {\bar\vt}

\nek{\kd} {\text{\large\tt Cod}(\id11)}

\nek{\bez} {\dif}

\nek{\fB} {\mathbf B}

\nek{\gh} {Ганди -- Харрингтон}
\nek{\bor} {{Borel}}

\nek{\Bore} [2] {{\prift Борелевск{#1} Продолжени{#2}}}

\nek{\bnd}{{(\bn)}{}^\dN}

\nek{\haf} [1] {\chi_{#1}}

\nek{\Gds} {\fG_{\fda\fsg}}
\nek{\Fsd} {\fF_{\fsg\fda}}

\nek{\prom} [1] {\text{\bfit P}(#1)}
\nek{\frg} [1] {\text{\sfit F\hspace*{0.05ex}}_{#1}}
\nek{\frd} [2] {(#1){}^{#2}}

\nek{\sfit}{\sffamily\itshape}

\nek{\goa} {{\BBB{\mathfrak a}\BBB}}
\nek{\gob} {{\BBB{\mathfrak b}\BBB}}

\nek{\np}{\newpage}

\nek{\br} {{\text{\bfit R}}}

\nek{\ls} {<_{\text{\sc лс}}}
\nek{\lse} {\le_{\text{\sc лс}}}

\nek{\ogr} [2] {{#1}\res_{#2}}

\nek{\ofu} [1] {{#1}_{\text{\tt wf}}}

\nek{\clo} [1] {\overline{#1}}

\nek{\pro} [2] {{#1}^{(#2)}}

\nek{\kbr} [1] {|#1|_{\text{\tt CB}}}

\nek{\hF} {\widehat F}
\nek{\ds} {\ddd{\fda s}}

\nek{\Siy} {\dS_\iy}

\nek{\raz} [1] {\widehat{#1}}

\nek{\bfm} {\mathbf m}

\nek{\num} {\text{\tt num}}

\nek{\Uni} [1] {{\prift Униформизаци{#1}}}

\renek{\xm} [1] {\dm#1\dm}
\renek{\qq} [1] {#1}

\nek{\fw}{\mathbf w}

\nek{\pop} [1] {\mathrel{\leqslant^*_{#1}}}
\nek{\pops} [1] {\mathrel{<^*_{#1}}}
\nek{\pom} [1] {\mathrel{\leqslant^-_{#1}}}
\nek{\poms} [1] {\mathrel{<^-_{#1}}}

\nek{\leqs} {\leqslant}

\nek{\lest} {\leqs^\ast}
\nek{\lst} {<^\ast}

\nek{\typ} [1] {|#1|}

\nek{\fk} {\boldsymbol k}

\nek{\cki} {\omega_1^{\text{\sc CK}}}
\nek{\ckp} [1] {\omega_1^{\text{\sc CK}(#1)}}

\nek{\cun} {\mathop\nabla}
\nek{\com} [1] {\neg\,{#1}}

\nek{\bi} {\mathbf 1}

\nek{\lcap} {\bigwedge}
\nek{\lcup} {\bigvee}

\nek{\rab} {\mathop{\tt rk}}

\nek{\cog} [1] {генерическ{#1} по Коэну}
\nek{\hcp} {\hc_{\gc^+}}
\nek{\tc} [1] {\mathrm{TC}(#1)}

\nek{\ua} {\dot a}
\nek{\uG} {\underline G}
\nek{\kf} [1] {{\mathbf C}_{#1}}
\nek{\sha} [2] {U_{#1}(#2)}
\nek{\kg} [1] {{\mathbf C}'_{#1}}
\nek{\shc} [2] {\overline U_{#1}(#2)}

\nek{\fA} {\mathbf{KA}}

\nek{\fV} {\mathbf{V}}
\nek{\kV} {\mathbb V}
\nek{\kL} {\mathbb L}
\nek{\gM} {{\mathfrak M}}
\nek{\bai} {\bn}
\nek{\gop} [1] {{\mathfrak F}_{#1}}
\nek{\dles}{\lessdot}

\nek{\el}  [2] {\xe_{#1}[#2]}
\nek{\El}  [2] {\xee_{#1}[#2]}
\nek{\eli} [1] {\xe_{#1}}
\nek{\Eli} [1] {\xee_{#1}}
\nek{\xee} {\boldsymbol{F}} %{{\mathfrak S}}
\nek{\xe} {\boldsymbol{f}}

\nek{\cor} [1] {\mathrel{{<}[#1]}}
\nek{\ubF}{\fontseries{b}\fontshape{it}\selectfont}
\nek{\xs} {\boldsymbol{s}}

\nek{\AC}   {{\text{\ubf AC}}}
\nek{\DC}   {{\text{\ubf DC}}}
\nek{\AD}   {{\text{\ubf AD}}}
\nek{\PD}   {{\text{\ubf PD}}}
\nek{\GCH}  {{\text{\ubf GCH}}}
\nek{\CH}   {{\text{\ubf CH}}}

\nek{\omj} [1] {\om_1^{\kL[#1]}}
\nek{\ko}  [1] {\mathrel{\prec_{#1}}}
\nek{\ok}  [1] {\mathrel{\succ_{#1}}}

\nek{\LR} [1]    {\kL[#1]\cap\bai}
\nek{\Lr}        {\kL\cap\bai}
\nek{\Or}        {\OD\cap\bai}

\nek{\enrm}
{\def\theenumi{{\rm(\roman{enumi})}}%
\itsep}

\nek{\smpr} {\hspace*{0.3ex}}

\nek{\bsg} {\smpr\bolf\sigmA\smpr}%{\text{\mtho\large$\cur\sigmA$}}
\nek{\bpi} {\smpr\bolf\pI\smpr}   %{\text{\mtho\large$\cur\pI$}}
\nek{\csg} {\bsg'}%{\text{\mtho\large$\cur\sigmA$}}
\nek{\cpi} {\bpi'}   %{\text{\mtho\large$\cur\pI$}}

\nek{\bfi} {\bolf\fI}   %{\text{\mtho\large$\cur\fI$}}
\nek{\bsi} {\bolf\psI}

\nek{\rve} {\mathrel{\varepsilon}}
\nek{\bolf} [1] {{\cur#1}} % {\text{\mtho\large$\cur#1$}}
\mathchardef\fI    ="7127
\mathchardef\psI   ="7120

\nek{\WO} {\wo}
\nek{\WF} [1] {{\text{\ubf WF}}_{#1}}
\nek{\pfr} [1] {п.\:\ref{#1}}
\nek{\mne} [1] {S(#1)} %{S^{{\rve}}_{#1}}
\nek{\bS} {{\bf S}}
\nek{\rss} [2] {#1\hspace*{0.1ex}{\res}\hspace*{0.2ex}_{#2}}
\nek{\BC} {\bok}

\nek{\atli} {\addtolength{\itemsep}{\dxi}}
\newlength{\dxi}
\nek{\axce}{\hspace{0\mathsurround}}
\nek{\axbx}[1]{\axce\protect\nolinebreak{\axfnt#1\/}%
\axce\protect\nolinebreak}
\nek{\hxbx}[1]{\axce\protect\nolinebreak{{\axfnt#1\/}\axce}%
\protect\nolinebreak}
\nek{\itea} [1] 
{\item[{\fontseries{m}\fontshape{sl}\selectfont#1\axce}%
{\fontshape{n}\selectfont:}]}

\nek{\orr} {\lor}  %mathbin{\tyl\bigvee}}

\nek{\Ext} [1] {\hxbx{Экс\-тен\-си\-о\-наль\-ност{#1}}}
\nek{\Pai} [1] {\hxbx{Пар{#1}}}
\nek{\Sep} [1] {\hxbx{Выделени{#1}}}
\nek{\Com} [1] {\hxbx{Объемност{#1}}}
\nek{\Rep} [1] {\hxbx{Подстановк{#1}}}
\nek{\Col} [1] {\axbx{Собирани{#1}}}
\nek{\Ifa} [1] {\hxbx{Бес\-ко\-неч\-ност{#1}}}
\nek{\Reg} [1] {\hxbx{Ре\-гу\-ляр\-ност{#1}}}
\nek{\Fou} [1] {\hxbx{Фундировани{#1}}}
\nek{\PS}  [1] {\axbx{Степен{#1}}}
\nek{\Cho} [1] {\hxbx{Выбор{#1}}}
\nek{\Wel} {\hxbx{Well-Ordering}}
\nek{\Uno} [1] {\hxbx{Объединени{#1}}}

\nek{\shas} {schemata}
\nek{\axfnt}{\fontfamily{cmss}\selectfont} 

\nek{\bp} [1] {\xbp(#1)}
\nek{\pk} [1] {\xpk(#1)}
\nek{\lm} [1] {\xlm(#1)}
\nek{\xlm} {{\text{\fontshape{n}\selectfont LM}}}
\nek{\xpk} {{\text{\fontshape{n}\selectfont PK}}}
\nek{\xbp} {{\text{\fontshape{n}\selectfont BP}}}

\nek{\kK} {{\bf K}}

\nek{\can} {\dn}
\nek{\fmu} {{\boldsymbol{\la}}}

\nek{\kP} {{\rm P}}
\nek{\bsp} {{\boldsymbol{p}}}

\nek{\Gen}[1]    {\text{\ubf Coh}\hspace*{0.35ex}{#1}}
\nek{\Ram}[1]    {\Ra\mu{#1}}
\nek{\Ra} [2]    {\text{\ubf Rand}_{#1}{\hspace*{0.2ex}#2}}
\nek{\Ran}[1]    {\Ra{\hspace*{0.2ex}}{#1}}
\nek{\Qa} [2]    {\Ra{#1}{\kL[#2]}}
\nek{\Qan}[1]    {\Qa{\hspace*{0.2ex}}{#1}}
\nek{\Qen}[1]    {\Gen{\kL[#1]}}

\nek{\enar} {\tenu{{\rm(\arabic{enumi})}}}
\nek{\enRm}
{\def\theenumi{(\Roman{enumi})}\itsep}
\nek{\enrmp}
{\def\theenumi
{{\mtho$(\text{\rm\roman{enumi}}\hspace{0.2ex}')$}}%
\itsep}

\nek{\tmu} [1] {\mathop{\tt cod}_{#1}} %{{\wh\mu}}

\nek{\mobox} [1] {\parbox[t]{\molen}{\mtho$#1$}}
\nek{\omk} {\om_1^{\kL}}

\nek{\tow} [1] {{тощ#1}}
\nek{\Tow} [1] {{Тощ#1}}
\nek{\kot} [1] {ко{тощ#1}}
\nek{\Kot} [1] {Ко{Тощ#1}}
\nek{\izam} {Исторические и библиографические замечания}

\nek{\kod} [1]  {\mathop{\text{\ubf cod}}#1}

\nek{\icat} {\cI_{\text{cat}}}

\nek{\Rai}[1]    {\Ra\cI{#1}}
\nek{\sgu} {\ddd{\sg\text{\rm-\linebreak[0]\ccc}}}

\nek{\enci} {\tenu{{\rm\mtho$\arabic{enumi}^\circ$}}
}

\nek{\se} [2] {{#1}^{(#2)}}

\nek{\vko} [1] {{\mtho$#1$-\yz в-\yz ко\-дах}}

\nek{\lja}[2] {\xxxm_{#1}(#2)}
\nek{\xxxm}  {{\text{\fontshape{n}\selectfont M}}}
\nek{\xgm} {{\text{\fontshape{n}\selectfont GM}}}

\nek{\dla} {\ddd{\kL[a]}}
\nek{\dlo} {\ddd{\kL}}
\nek{\lia}[1] {\xgm(#1)}

\nek{\ency} {\tenu{{\rm(\asbuk{enumi})}}}
\nek{\encyp}
{\def\theenumi
{{\mtho$(\text{\rm\asbuk{enumi}}\hspace{0.2ex}')$}}%
\itsep}

\nek{\pkc}[1] {\xpk^-(#1)}

\nek{\iza} [1]
{\subsubsection*{{\bfit\izam}}}

\nek{\raj}[1]    {\Ra\cI{\kL[#1]}}

\nek{\gfr} {\nrf}
\nek{\ROD}{{\text{\ubf ROD}}}
\nek{\hrod}{{\text{\ubf HROD}}}

\nek{\sekd}{\punk}
\nek{\usl} [1] {\lap{услови#1}}
\nek{\Usl} [1] {\lap{Услови#1}}

\nek{\enaR} {\tenu{{\rm\arabic{enumi})}}}

\nek{\ter} [1] {{\breve{#1}}}
\nek{\tex} {{\ter x}}
\nek{\tey} {{\ter y}}
\nek{\tez} {{\ter z}}

\nek{\vy} {\mathrel{|\hspace*{-0.1ex}|\hspace*{-0.5ex}{-}%
\hspace*{-1.0ex}{-}}}

\nek{\vyn} [2] {\vy^{#1}_{#2}}
\nek{\vin} [1] {\vyn{}{#1}}

\nek{\dCx} {\dC({\xi})}
\nek{\pli} {\mathord{\hspace*{0.1ex}+\hspace*{0.1ex}}}
\nek{\mii} {\mathord{\hspace*{0.1ex}-\hspace*{0.1ex}}}

\nek{\Soi} {{\Sol\cI}}
\nek{\soi} [1] {\sol{#1}\cI}
\nek{\sol} [2] {\dP^{#1}_{#2}}
\nek{\qol} [2] {\sol{\kL[#1]}{#2}}
\nek{\Sol} [1] {\dP_{#1}}

\nek{\gel} [1] {\pi_{#1}}
\nek{\geg} {{\gel G}}

\nek{\lpo} [1] {{\mathord{^{#1}\hspace*{-0.2ex}2}}}
\nek{\qer} [1] {{\breve{\phantom{x}}\hspace*{-0.5ex}{#1}}}

\nek{\upi} {{\underline\pi}}
\nek{\Coh}     {\dP_{\text{\tt coh}}}
\nek{\coh} [1] {\dP_{\text{\tt coh}}^{#1}}
\nek{\Som}     {\Sol\fmu}
\nek{\som} [1] {\sol{#1}\fmu}

\nek{\yC}     {{\mathbf C}}
\nek{\yD}     {{\mathbf D}}

\nek{\sqi} {\sq_\cI}
\nek{\eqi} {\approx}
\nek{\para} [2] {{#1}{\hspace*{0.2ex}\tilde{\phantom{.}}}{#2}}
\nek{\lapr}[1] {{\rm\lap{#1}}}
\nek{\qoi} [1] {\qol{#1}\cI}

\nek{\Tr} {{\ubf{Tr}}}

\nek{\lee} {\le^\ast}
\nek{\ine} {\in^\ast}

\nek{\fell} {\text{\usefont{OML}{cmm}{b}{it}\char"60}}
\nek{\feli} {{\fell_1}}
\nek{\qai} {{{\dQ_+}^\om}}
\nek{\kel} [1] {\##1}
\nek{\lsup} {\mathop{\text{\tt lim\,sup}}}
\nek{\ki} [1] {\bon{#1}\dn}
\nek{\Bi} [1] {\bon{#1}\bn}

\nek{\OD} {{\text{\rm OD}}}
\nek{\ZFI}  {{\text{\ubf ZFCI}}}
\nek{\od}   {\mathop{\text{\tt od}}}
\nek{\rod}  {\mathop{\text{\tt rod}}}

\nek{\dCi} {\dC({\omk})}
\nek{\dCt} {\dC({\vt})}
\nek{\dCk} {\dC({\ka})}
\nek{\tea} {{\ter a}}
\nek{\tek} {{\ter k}}
\nek{\ten} {{\ter n}}
\nek{\bP} {{\mathbf{P}}}

\nek{\dagg} {{\mtho(\ensuremath\dag)}}
\nek{\dagd} {{\mtho(\ensuremath\ddag)}}
\nek{\sP} {{\cP}}

\nek{\Coll}  {{\text{\rm Coll}}}
\nek{\hgt}  [1]  {\mathop{\text{\tt hgt}}#1}

\nek{\ana} [1] {\boldsymbol A[#1]}
\nek{\coa} [1] {\boldsymbol C[#1]}
\nek{\ank} [2] {\boldsymbol A_{#2}[#1]}
\nek{\cok} [2] {\boldsymbol C_{#2}[#1]}

\nek{\yA} [1] {{\mathbf A}_{#1}}
\nek{\yB}     {{\mathbf B}}
\nek{\zA} [2] {(\yA{#1})^{\kL[#2]}}
\nek{\zB} [1] {(\yB)^{\kL[#1]}}
\nek{\zC}     {\yC}
\nek{\zD} [1] {(\yD)^{\kL[#1]}}

\nek{\MA}  {\text{\ubf MA}}

\nek{\Ts}  {{\overline s}}
\nek{\tQ} {{\overline\dQ}}

\nek{\iyf} {{\imy f}}
\nek{\imy} {\poq}

\nek{\teg} {{\ter g}}
\nek{\tej} {{\ter j}}

\nek{\kan} [1] {\cC^{#1}}

\nek{\jx} {{\boldsymbol{x}}}
\nek{\jy} {{\boldsymbol{y}}}
\nek{\ja} {{\boldsymbol{a}}}

\nek{\bss}[2] {{\cur B}_{#2}^{#1}}
\nek{\bsl}[1] {\bss{\vt}{#1}}
\nek{\bsx}[1] {\bss{\xi}{#1}}

\nek{\bea} [1] {\cN^{#1}}

\nek{\gds} {\mathbf U}

\nek{\jP} {{\boldsymbol{P}}}

\nek{\bod} [2] {U^{#1}_{#2}}

\nek{\cjj} [1] {\cI_{#1}}

\nek{\ctt} {{\mathbf t}}

\nek{\bc} {\text{\ubf bc}}
\nek{\kos} {\fsg}
\nek{\kop} {\fpi}

\nek{\rL}{\text{\rm L}}
\nek{\sqq}{\sqsubseteq}

\nek{\oE}  {\mathbin{\overline{\rE}}}
\nek{\noE} {\mathbin{\not{\hspace{-2pt}\overline{\rE}}}}
\nek{\nEo} {\mathbin{\not{\hspace{-2pt}\Eo}}}

\nek{\pri}  {{\tt pr}_1\hspace{1.5pt}}
\nek{\prii} {{\tt pr}_2\hspace{1.5pt}}

\nek{\J} [1] {\mathbin{\rR_{#1}}}
\nek{\I} [1] {\mathbin{\rQ_{#1}}}
\nek{\Ip}[1] {\mathbin{\rQ'_{#1}}}

\nek{\cXp} {\cX^\ast}
\nek{\cPp} {\cP^\ast}
\nek{\mins}{\hspace{-1pt}-\hspace{-1pt}}

\nek{\oet} {{{\oE}{\hspace{0.1ex}}^3}} 
\nek{\oed} {\oE} 
\nek{\lxx} {<^{\ast\ast}} 
\nek{\lx}  {<^\ast} 
\nek{\mex} {\leqslant^\ast} 
\nek{\gex} {\geqslant^\ast} 
\nek{\ex}  {\equiv^\ast}

\nek{\gap}[3]{\ddd{(#1,#2)}{#3}}
\nek{\kpa} {\ka}

\nek{\hau} [2] {(#1\:#2)} 

\nek{\gam} [1] {{\mathbf G}_{#1}}
\nek{\gan} [3] {{\mathbf G}_{#1}(#2\hspace*{0.2ex};\hspace*{0.1ex}#3)}

\nek{\pos} [2] {\ang{#1\hspace{0.2ex};\linebreak[0]\hspace{0.1ex}#2}}

\nek{\gq} {\mathop{{\Game}{\hspace*{-1.44ex}}{\Game}%
\hspace*{-1.44ex}{\Game}\hspace*{-1.44ex}{\Game}}}

\nek{\<} {\le}

\mathchardef\pI ="7119
\theoremstyle{plain}

\newtheorem{princip}     [theore]{Принцип}
\newtheorem{fact}        [theore]{Факт}

\nek{\bpri}{\begin{princip}}
\nek{\epri}{\end{princip}}
\nek{\bfa}{\begin{fact}}
\nek{\efa}{\end{fact}}

\nek{\as}[1]{{\displaystyle\iSg^\dag_{#1}}}
\nek{\ar}[1]{{\displaystyle\iPi^\dag_{#1}}}

\nek{\ms}[1]{{\displaystyle\iSg^\ast_{#1}}}
\nek{\mr}[1]{{\displaystyle\iPi^\ast_{#1}}}

\DeclareFontFamily{U}{Bcmm}{\skewchar\font127 }
\DeclareFontShape{U}{Bcmm}{m}{n}{ <->[1.315]cmmi10}{}%
\DeclareFontShape{U}{Bcmm}{b}{n}{ <->[1.315]cmmib10}{}%
\DeclareMathAlphabet{\nmi}{U}{Bcmm}{m}{n}%
\DeclareMathAlphabet{\bmi}{U}{Bcmm}{b}{n}%

\nek{\kp} [1] {{\bmi\pI}_{#1}} 
\nek{\jp} [1] {{\nmi\pI}_{#1}} 

\nek{\bs} [1] {\boldsymbol B_{#1}}
\nek{\bt} [2] {\boldsymbol B_{#1}(#2)}

\nek{\wX} {\widehat X}

\nek{\bok} {\text{\ubf BC}} 

\nek{\zfct} {\ZFC^-_{\alo}}
\nek{\zhc}  {\zfct} 

\nek{\op} {{\boldsymbol p}}
\nek{\ou} {{\boldsymbol u}}
\nek{\ov} {{\boldsymbol v}}
\nek{\oa} {{\boldsymbol a}}
\nek{\ob} {{\boldsymbol b}}
\nek{\ox} {{\boldsymbol x}}

\nek{\de} [1] {\text{\sc Det}(#1)} 

\nek{\bnt} {(\bn){}^2}

\nek{\BM} [1] {\text{\rm BM}(#1)}

\nek{\stra} {{\skr S}}

\nek{\fT} {\boldsymbol T}

\nek{\igri} [4] {G^{#4}_{#1}(#2,#3)}
\nek{\igrd} [5] {G^{#4,#5}_{#1}(#2,#3)}

\nek{\pra} [2] {{\prec}#1,#2{\succ}}

\nek{\leuv} [3] {\le_{#1#2#3}}

\nek{\bis} [2] {\text{\ubf BS}_{#1\hspace{0.1ex}#2}}
\nek{\bisi} {\text{\ubf BiSi}}

\nek{\gom} {\text{\ubf Hom}}
\nek{\gor} {\text{\ubf Hom}^{\text{\rm reg}}}

\nek{\pho} [2] {\text{\ubf PH}_{#1#2}}
\nek{\rhom} [2] {\text{\ubf RH}_{#1#2}}
\nek{\rpho} [2] {\text{\ubf PH}^{\text{\rm reg}}_{#1#2}}

\nek{\aop}   {A-операци}
\nek{\ase}   {A-множеств}
\nek{\cse}   {CA-множеств}
\nek{\biz}   {B-измерим}

\nek{\mon} {\boldsymbol n}
\nek{\mom} {\boldsymbol m}
\nek{\balo} {\boldsymbol\aleph_0}
\nek{\bcon} {\boldsymbol\cont}

\nek{\qw} {\mathopen{\text{\sf W}\hspace{0.25ex}}}

\nek{\bX} {\mathbf X}

\nek{\puns} {\subsubsection}
\renek{\thesubsubsection}
{\arabic{subsubsection}}

\renek{\punk} {\puns}

\makeatletter
\if@twoside%
\else\fi%
\makeatother  

\begin{document}

\title{Linear ROD subsets of Borel partial orders are 
countably cofinal in Solovay's model}

\author{Vladimir Kanovei~\thanks
{IPPI, Moscow, Russia.}
}

\date{\today}
\maketitle

\begin{abstract}
The following is true in the Solovay model.
  
1. If $\stk D\le$ is a Borel partial order on a set $D$ of the 
reals, $X\sq D$ is a ROD set, and ${\le}\res X$ is 
linear, then ${\le}\res X$ is countably cofinal.  

2. If in addition every countable set $Y\sq D$ has a strict upper 
bound in $\stk D\le$ then the ordering $\stk D\le$ has no maximal 
chains that are ROD sets.
\end{abstract}

Linear orders, 
which typically appear in conventional mathematics,  
are countably cofinal. 
In fact \rit{any} Borel (as a set of pairs) linear order on 
a subset of a Polish space is countably cofinal: 
see, \eg, \cite{hms}. 
On the other hand, there is an uncountably-cofinal quasi-order of  
class $\fs11$ on $\bn$.

\bex
\lam{ex}
Fix any recursive enumeration $\dQ=\ens{q_k}{k\in\dN}$ 
of the rationals. 
For any ordinal $\xi<\omi$, let $X_\xi$ 
be the set of all points $x\in\bn$ such that the maximal well-ordered 
(in the sense of the usual order of the rationals) 
initial segment of the set 
$
Q_x=\ens{q_k}{x(k)=0}
$ 
has the order type $\xi$. 
Thus $\bn=\bigcup_{\xi<\omi}X_\xi$. 
For $x,y\in\bn$ define $x\lef y$ iff $x\in X_\xi$, 
$y\in X_\eta$, and $\xi\le \eta$. 
Thus $\lef$ is a prewellordering of length exactly $\omi$. 
It is a routine exercise to check that $\lef$ belongs to $\fs11$. 

We can even slightly change the definition of $\lef$ to obtain 
a true linear order.  
Define $x\lef' y$ iff either $x\in X_\xi$, 
$y\in X_\eta$, and $\xi<\eta$, or $x,y\in X_\xi$ for one 
and the same $\xi$ and $x<y$ in the sense of the lexicographical 
linear order on $\bn$. 
Clearly $\lef'$ is a linear order of cofinality $\omi$ and 
class $\fs11$. 
\eex

Yet there is a rather representative 
class of \ROD\ (that is, real-ordinal definable) 
linear orderings which are 
consistently countably cofinal. 
This is the subject of the next theorem.\pagebreak[4]

\bte
\lam{m}
The following sentence is true in the Solovay model$:$
if\/ $\le$ is a Borel partial quasi-order on a (Borel) set\/ 
$D\sq\bn$, $X\sq D$ is a\/ \ROD\ set, and\/ 
${\le}\res X$ is a linear quasi-order, then\/ ${\le}\res X$ is 
countably cofinal.
\ete

A \rit{partial quasi-order}, PQO for brevity, is a 
binary relation $\le$ satisfying 
${x\le y}\land {y\le z}\imp {x\le z}$ and $x\le x$ 
on its domain. 
In this case, an \rit{associated equivalence relation} $\equiv$ 
and an \rit{associated strict partial order} $<$ are defined 
so that $x\equiv y$ iff ${x\le y}\land {y\le x}$, and 
$x<y$ iff ${x\le y}\land {y\not\le x}$.
A PQO is \rit{linear}, LQO for brevity, if we have 
${x\le y}\lor {y\le x}$ for all $x,y$ in its domain.

A PQO $\stk X\le$ 
(meaning: $X$ is the domain of $\le$) 
is \rit{Borel} iff the set $X$ is 
a Borel set in a suitable Polish space $\dX$, and 
the relation $\le$ 
(as a set of pairs) is a Borel subset of $\dX\ti\dX$.

Thus it is consistent with $\ZFC$ that \ROD\ linear suborders 
of Borel PQOs are necessarily countably cofinal. 
Accordingly it is consistent with $\ZF+\DC$ that any linear suborders 
of Borel PQOs are countably cofinal.

By \rit{the Solovay model} we understand a model of $\ZFC$ in 
which all \ROD\ sets of reals have some 
basic regularity properties, for instance, are Lebesgue measurable, 
have the Baire property, see \cite{sol}.
We'll make use of the following two results related to the 
Solovay model.

\bpro
[Stern~\cite{stl}]
\lam{p2}
It holds in the Solovay model that if\/ $\rho<\omi$ 
then there is no\/ \ROD\/ \ddd{\omi}sequence of pairwise different   
sets in\/ $\fs0\rho$. 
\qed  
\epro

\bpro
%[Kanovei~\cite{aord}]
\lam{p3}
It holds in the Solovay model that if\/ $\le$ is 
a\/ \ROD\ LQO\ 
on a set\/ $D\sq\bn$ then there exist a\/ \ROD\ antichain\/ 
$A\sq2^{<\omi}$ and a\/ \ROD\ map\/ $\vt:D\onto A$  
such that\/ $x\le y\eqv \vt(x)\lexe\vt(y)$ for all\/ $x,y\in D$.
\qed
\epro

A few words on the notation. 
The set $2^{<\omi}=\bigcup_{\xi<\omi}2^\xi$ 
consists of all transfinite binary sequences of length $<\omi$, 
and if $\xi<\omi$ then $2^\xi$ is the set 
of all binary sequences of length exactly $\xi$. 
A set $A\sq2^{<\omi}$ is an \rit{antichain} if we have $s\not\su t$ 
for any $s,t\in A$, where $s\su t$ means that $t$ is a proper 
extension of $s$. 
By $\lexe$ we denote the lexicographical order on $2^{<\omi}$, 
that is, if $s,t\in2^{<\omi}$ then $s\lexe t$ iff 
either 1) $s=t$ or 
2) $s\not\su t$, $t\not\su s$, and the least ordinal 
$\xi<\dom s\yi \dom t$ such that $s(\xi)\ne t(\xi)$ satisfies 
$s(\xi)<t(\xi)$.
Obviously $\lexe$ linearly orders any antichain $A\sq 2^{<\omi}$.

Proposition~\ref{p3} follows from Theorem 6 in \cite{aord} saying that 
if, in the Solovay model, $\le$ is a\/ \ROD\ PQO\ on a set\/ $D\sq\bn$
then:
\bde
\item[\it either] 
a condition $({\rm I}^s)$ holds, which for LQO relations $\le$  
is equivalent to the existence of $A$ and $\vt$ as in Proposition~\ref{p3}, 

\item[\it or] 
a condition $({\rm II})$ holds, which is incompatible with 
$\le$ being a LQO.
\ede
Thus we obtain Proposition~\ref{p3} as an immediate corollary.

The next simple fact will be used below.

\ble
\lam{cc}
If\/ $\xi<\omi$ then any set\/ $C\sq2^\xi$ is countably\/ 
\ddd\lexe cofinal, 
that is, there is a set\/ $C'\sq C$, at most countable and\/ 
\ddd\lexe cofinal in\/ $C$.\qed
\ele 

\bpf[Theorem~\ref{m}]
\rit{We argue in the Solovay model.}
Suppose that $\le$ is a Borel PQO on a (Borel) set $D\sq\bn$, 
$X\sq D$ is a \ROD\ set, and ${\le}\res X$ is a LQO. 
Our goal will be to show that ${\le}\res X$ is 
countably cofinal, that is, there is a set 
$Y\sq X$, at most countable and \ddd{\le}cofinal in $X$.

The restricted order ${\le}\res X$ is \ROD, of course, and 
hence, by Proposition~\ref{p3}, 
there is a \ROD\ map $\vt:X\onto A$ onto an antichain $A\sq2^{<\omi}$ 
(also obviously a \ROD\ set)
such that $x\le y\eqv \vt(x)\lexe\vt(y)$ for all $x,y\in X$. 

If $\xi<\omi$ then let $A_\xi=A\cap 2^\xi$ 
and $X_\xi=\ens{x\in D}{\vt(x)\in A_\xi}$.\vom

\rit{Case 1}: there is an ordinal $\xi_0<\omi$ such that 
$A_{\xi_0}$ is \ddd\lexe cofinal in $A$.
However, by Lemma~\ref{cc}, there is a set 
$A'\sq A_{\xi_0}$, at most countable and \ddd\lexe cofinal in 
$A_{\xi_0}$, and hence \ddd\lexe cofinal 
in $A$ as well by the choice of $\xi_0$. 
If $s\in A'$ then pick an element $x_s\in X$ such that $\vt(x_s)=s$. 
Then the set $Y=\ens{x_s}{s\in A'}$ is a countable subset of $X,$ 
\ddd\le cofinal in $X,$ 
as required.\vom

\rit{Case 2}: not Case 1. 
That is, for any $\eta<\omi$ there is an ordinal $\xi<\omi$ and 
an element $s\in A_\xi$ such that $\eta<\xi$ and $t\lex s$ 
for all $t\in A_{\eta}$.
Then the sequence of sets 
$$
D_\xi=\ens{z\in D}{\sus x\in X\,(z\le x\land \vt(x)\in A_\xi)}
$$
%$B_\xi=\bigcup_{\eta<\xi}A_\xi$ 
is \ROD\ and has uncountably many pairwise different terms. 
%by the Case~2 assumption. 

We are going to get a contradiction.
Recall that $\le$ is a Borel relation, hence it belongs to 
$\fs0\rho$ for an ordinal $1\le\rho<\omi$.
Now the goal is to prove that all sets $D_\xi$ belong to 
$\fs0\rho$ as well --- this contradicts to Proposition~\ref{p2}, 
and the contradiction accomplishes the proof of the theorem. 

Consider an arbitrary ordinal $\xi<\omi$. 
By Lemma~\ref{cc} there exists a countable set 
$A'=\ens{s_n}{n<\om}\sq A_\xi$, 
\ddd\lexe cofinal in $A_\xi$.
If $n<\om$ then pick an element $x_n\in X$ such that 
$\vt(x_n)=s_n$. 
Note that by the choice of $\vt$ any other element $x\in X$ 
with $\vt(x)=s_n$ satisfies $x\equiv x_n$, where $\equiv$ 
is the equivalence relation on $D$ associated with $\le$. 
It follows that
$$
\textstyle
D_\xi=\bigcup_nX_n\,,\quad\text{where}\quad
X_n=\ens{z\in D}{z\le x_n}\,,
$$
so each $X_n$ is a $\fs0\rho$ set together with $\le$, and so 
is $D_\xi$ as a countable union of sets in $\fs0\rho$.\vom

\epF{Theorem~\ref{m}}\vom

We continue with a few remarks and questions. 

\bqe
Can one strengthen Theorem~\ref m as follows: 
\rit{the restricted relation\/ 
${\le}\res X$ has no monotone\/ \ddd\omi sequences}?
Lemma~\ref{cc} admits such a strengthening: 
if $\xi<\omi$ then easily any \ddd\lexe monotone 
sequence in $2^\xi$ is countable.
\eqe

Using Shoenfield's absoluteness, we obtain:

\bcor
\lam{mc}
If\/ $\le$ is a Borel PQO on a (Borel) set\/ 
$D\sq\bn$, $X\sq D$ is a\/ $\fs11$ set, and\/ 
${\le}\res X$ is a linear quasi-order, 
then\/ ${\le}\res X$ is 
countably cofinal.
\ecor

Note that Corollary~\ref{mc} fails for arbitrary LQOs of class $\fs11$ 
(that is, not necessarily linear suborders of Borel PQOs), see 
Example~\ref{ex}.

\bpf
In the case considered, the property of countable cofinality of 
${\le}\res X$ can be expressed by a $\is12$ formula. 
Thus it remains to consider a Solovay-type extension of the 
universe and refer to Theorem~\ref m.\snos 
{We'll not discuss the issue of 
an inaccessible cardinal on the background.}
\epf

Yet there is a really elementary proof of Corollary~\ref{mc}. 

Let $Y$ be the set of all elements $y\in D$ \ddd\le comparable with 
\rit{every} element $x\in X$. 
This is a $\fs11$ set, and $X\sq Y$ 
(as ${\le}$ is linear on $X$). 
Therefore there is a Borel set $Z$ such that $X\sq Z\sq Y$. 
Now let $U$ be the set of all $z\in Z$ \ddd\le comparable with 
\rit{every} element $y\in Y$. 
Still this is a $\fs11$ set, and $X\sq U$ by the definition of $Y$.
Therefore there is a Borel set $W$ such that $X\sq W\sq U$. 
And by definition still ${\le}$ is linear on $W$. 
It follows that $W$ does not have increasing \ddd\omi sequences, 
and hence neither does $X$.

\bqe
Is Corollary~\ref{mc} true for $\fp11$ sets $X$?

We cannot go much higher though. 
Indeed, if $\le$ is, say, the eventual domination order on $\bn$, 
then the axiom of constructibility implies the existence of a 
\ddd\le monotone \ddd\omi sequence of class $\id12$. 
\eqe

Now a few words on Borel PQOs $\le$ having the 
following property:
\ben
\fenu
\itla{*}
if $X$ is a countable set in the domain of $\le$ then there is an 
element $y$ such that $x<y$ 
(in the sense of the corresponding strict ordering)
for all $x\in X$. 
\een
A thoroughful study of some orderings of this type 
(for instance, the ordering on $\dR^\om$ defined so that 
$x\le y$ iff either $x(n)=y(n)$ for all but finite $n$ or 
$x(n)<y(n)$ for all but finite $n$) 
was undertaken in early papers of Felix Hausdorff, \eg, \cite{h07,h09} 
(translated to English in \cite{h05}). 
In particular, Hausdorff investigated the structure of 
\rit{pantachies}, that is, maximal 
linearly ordered subsets of those partial orderings. 
As one of the first explicit applications of the axiom of choice, 
Hausdorff established the existence of a pantachy in any partial 
order, and made clear distinction between such an existence proof 
and an actual, well-defined construction of an individual pantachy 
(see \cite{h07}, p.~110).
The next result shows that the latter is hardly possible in \ZFC, at least 
if we take for granted that any individual set-theoretic construction 
results in a \ROD\ set.

\bcor
\lam{l}
The following sentence is true in the Solovay model$:$
if\/ $\le$ is a Borel partial quasi-order on a (Borel) set\/ 
$D\sq\bn$, satisfying\/ \ref{*}, then\/ $\le$ has no\/ \ROD\ 
pantachies.
\ecor
\bpf
It follows from \ref{*} that any pantachy in $\stk D\le$ 
is a set of uncountable cofinality. 
Now apply Theorem~\ref m.
\epf

A further corollary: it is impossible to prove the existence of 
pantachies in any Borel PQO satisfying \ref{*} in {\ZF} + \DC.

\def\refname

\end{document}

\vyk{
3. 
We don't know whether Corollary~\ref{mc} is true for arbitrary LQOs 
of class $\fs11$. 
The following is a near-counterexample.
Let $D$ be the set of all trees $T\sq2\lom$. 
If $S,T\in D$ then define $S\le T$ iff there is a homomorphism 
$h:S\to T$. 
(That is, we require that $s\su t\imp h(s)\su h(t)$ for all strings 
$s\yi t\in S$, where $s\su t$ 
means that $t$ is a proper extension of $s$.)
Then $D$ is Borel, $\le$ is $\fs11$, and $S\le T$ iff 
either 
$T$ is ill-founded (and $S$ arbitrary) 
or 
both $S$ and $T$ are well-founded and the rank of $S$ is less or equal 
to the rank of $T$. 
Thus $\le$ is a prewellordering of length $\omi+1$ 
(with the highest class equal to the set of all ill-founded trees) 
and of class $\fs11$. 
It is still countably-cofinal, yet it admits an increasing 
\ddd\omi sequence.\vom  

4.
This leads to the question: 
is there a prewellordering $\le$ of class $\fs11$ and of length 
exactly $\omi$? 
Or somewhat weaker: 
is there a LQO $\le$ of class $\fs11$ and of cofinality  
exactly $\omi$? 
}

%\newpage